\documentclass[smallextended,envcountsect]{svjour3} 
\usepackage{latexsym, amsmath, amsfonts, amscd, amssymb, verbatim}
\usepackage{graphicx,graphics}
\usepackage{xcolor}
\smartqed 

\newtheorem{Theorem}{Theorem}[section]
\newtheorem{Proposition}{Proposition}[section]
\newtheorem{Remark}{Remark}[section]
\newtheorem{Lemma}{Lemma}[section]
\newtheorem{Example}{Example}[section]
\def\R{\rm I\!R}

\def\gph{\mbox{\rm gph}\,}

\def\Limsup{\mathop{{\rm Lim}\,{\rm sup}}}

\def\tto{\;{\lower 1pt \hbox{$\rightarrow$}}\kern -10pt
	\hbox{\raise 2pt \hbox{$\rightarrow$}}\;}

\def\R{I\!\!R}

\def\ox{\bar{x}}
\def\oy{\bar{y}}
\def\oz{\bar{z}}
\def\ov{\bar{v}}

\def\gph{\mbox{\rm gph}\,}

\begin{document}

\title{Sensitivity Analysis of a Stationary Point Set Map under Total Perturbations. Part 2: Robinson Stability}


\author{D. T. K.~Huyen  \and  J.-C.~Yao \and N.~D.~Yen}

\institute{D. T. K. Huyen \at Graduate Training Center, Institute of Mathematics, Vietnam Academy of Science and Technology,\\ Hanoi, Vietnam\\
 kimhuyenhy@gmail.com
 \and 
 Jen-Chih Yao \at Center for General Education, China Medical University,\\ Taichung, Taiwan\\
 yaojc@mail.cmu.edu.tw
 \and
	Nguyen Dong Yen, Corresponding author \at Institute of Mathematics, Vietnam Academy of Science and Technology,\\ Hanoi, Vietnam\\ ndyen@math.ac.vn}

\date{{\text{Communicated by ...}}
		\\
		\\
Received: date / Accepted: date}

\maketitle


\smallskip
\begin{abstract}
In Part 1 of this paper, we have estimated the Fr\'echet coderivative and the Mordukhovich coderivative of the stationary point set map of a smooth parametric optimization problem with one smooth functional constraint under total perturbations. From these estimates, necessary and sufficient conditions for the local Lipschitz-like property of the map have been obtained. In this part,  we establish sufficient conditions for the Robinson stability of the stationary point set map. This allows us to revisit and extend several stability theorems in indefinite quadratic programming. A comparison of our results with the ones which can be obtained via another approach is also given. 
\end{abstract}



\setcounter{equation}{0}
\section{Introduction}
Appeared at the early stage of optimization theory, smooth programming problems continue to attract common attention of the optimization community due to their importance and beauty. Polynomial optimization problems, including nonconvex quadratic programs, are typical examples of such problems.

\smallskip
The present paper investigates the Lipschitz-like property and the Robinson stability of the stationary point set map of a smooth parametric optimization problem with one smooth functional constraint under total perturbations.

\smallskip
In Part 1 of the paper \cite{Huyen_Yao_Yen_2018a}, we have computed and estimated the Fr\'echet coderivative and the Mordukhovich coderivative of the stationary point set map by applying some theorems of Levy and Mordukhovich \cite{LeMo04} and other related results. From the obtained formulas we derive necessary and sufficient conditions for the local Lipschitz-like property of the stationary point set map. This leads us to new insights into the preceding deep investigations of Levy and Mordukhovich in the just-cited paper and of Qui \cite{Qui_JOTA2014,Qui_JoGO_2016}. 

\smallskip
The reader is referred to Part 1 of this paper \cite{Huyen_Yao_Yen_2018a} for a survey on the local Lipschitz-like property of multifunctions, the Robinson stability of an implicit multifunction, the Mordukhovich criterion for the local Lipschitz-like property of locally closed multifunctions, and some relevant material.

\smallskip 
This part of the paper is organized as follows.  Section 2 recalls some basic concepts from variational analysis, formulates the problem studied herein, and presents a series of auxiliary results in a unified form. In Section 3, we obtain sufficient conditions for the Robinson stability of the stationary point set map. Section 4 is devoted to several stability theorems in indefinite quadratic programming. A comparison of our results with the ones which can be obtained via Robinson's theory of strongly regular generalized equations \cite{Robinson_1980} is given in Section 5. The final section contains some concluding remarks.

\section{Preliminaries}

The scalar product and the norm in a finite-dimensional Euclidean space are denoted respectively by $\langle
\cdot,\cdot \rangle$ and $\|\cdot \|$. The symbols $B(x,\rho)$ and $\bar{B}(x,\rho)$ stand for the open (resp., closed) ball centered at $x \in X$ with radius $\rho > 0$. The distance $\displaystyle\inf_{u\in A}\|x-u\|$ from $x \in X$ to a subset $A \subset X$ is denoted by $d(x, A)$.

We now recall several basic concepts from variational analysis \cite{B-M06,Rock_Wets_1998} which will be used intensively later on.

The {\it Fr\'echet normal cone} (also called the {\it prenormal
	cone}, or the {\it regular normal cone}) to a set
$\Omega\subset\R^s$ at $\ov \in \Omega$ is given by
\begin{eqnarray*}\label{Frechet normals} \widehat
	N_{\Omega}(\ov)=\left\{v'\in\R^s\, \mid\,
	\displaystyle\limsup_{v\xrightarrow{\Omega}\ov}\,\displaystyle\frac{\langle
		v',v-\ov \rangle}{\|v-\ov \|}\leq 0\right\},\end{eqnarray*} where
$v\xrightarrow{\Omega}\ov$ means $v\to \ov$ with $v\in\Omega$. By
convention, $\widehat N_{\Omega}(\ov):=\emptyset$ when $\ov \notin \Omega$.
Provided that $\Omega$ is locally closed around $\bar v\in \Omega$, one calls
\begin{eqnarray*}\label{basic normals}\begin{array}{rl}
		N_{\Omega}(\ov)&=\displaystyle\Limsup_{v\to\ov}\widehat
		N_{\Omega}(v)\\
		&:=\big\{v'\in\R^s\, \mid\,  \exists \mbox{ sequences } v_k\to \ov,\ v_k'\rightarrow v',\\
		& \qquad \qquad \qquad \quad \ \mbox {with } v_k'\in \widehat
		N_{\Omega}(v_k)\, \mbox{ for all }\, k=1,2,\dots
		\big\}\end{array}\end{eqnarray*} the {\it Mordukhovich} (or {\it
	limiting/basic}) {\it normal cone} to $\Omega$ at $\ov$.
If $\ov \notin \Omega$, then one puts $N_{\Omega}(\ov)=\emptyset$.

A multifunction $\Phi:\R^n\rightrightarrows  \R^m$ is
said to be {\it locally closed} around a point $\oz=(\ox,\oy)$ from  $\gph\Phi:=\{(x,y)\in \R^n\times\R^m\,\mid\, y\in\Phi(x)\}$ if $\gph\Phi$ is locally closed around $\bar z$. Here, the product space
$\R^{n+m}=\R^n\times\R^m$ is equipped with the topology generated by the sum norm $\|(x,y)\| =\|x\|+\|y\|$.

For any $\bar{z}=(\bar{x},\bar{y})\in \mbox{gph}\,\Phi$,
\begin{eqnarray*}\label{Frechet coderivative}
	\widehat D^*\Phi(\bar{z})(y'):=\big\{x'\in
	\R^n\, \mid\, (x',-y')\in \widehat N_{\mbox{gph}\,\Phi}(\bar{z})\big\}\quad (y'\in \R^m)
\end{eqnarray*} are called the {\it
Fr\' echet coderivative} values of $\Phi$ at $\bar z$. Similarly, the {\it
Mordukhovich coderivative} (limiting coderivative) values of $\Phi$ at $\bar z$ are defined by
\begin{eqnarray*}\label{normal coderivative}D^*\Phi(\bar{z})(y'):=\big\{x'\in
	\R^n\, \mid\, (x',-y')\in N_{\mbox{gph}\,\Phi}(\bar{z})\big\}\quad (y'\in \R^m).\end{eqnarray*}
Thus, $\widehat D^*\Phi(\bar{z})$ and $D^*\Phi(\bar{z})$ are multifuntions from $\R^m$ to $\R^n$. By \cite[Theorem~1.38]{B-M06}, if $\Phi$ is strictly Fr\' echet differentiable at $\bar x$, then
$$\widehat D^*\Phi(\bar{x})(y')=D^*\Phi(\bar{x})(y')
=\{\nabla\Phi(\bar x)^*(y')\}$$
for any $y'\in \R^m$.

Suppose that $X$, $Y$, and $Z$ are finite-dimensional Euclidean spaces. Consider a function $\psi: X\rightarrow \bar{\R}$ with $|\psi(\bar x)|<\infty$. The set
$$\partial\psi(\bar x):=\{x'\in X^*\mid (x',-1)\in N_{{\rm epi}\,\psi}(\bar x,\psi(\bar x))\}$$
is the \textit{Mordukhovich subdifferential} of $\psi$ at $\bar x$. We put $\partial\psi(\bar x)=\emptyset$ if $|\psi(\bar x)|=\infty$. The set
$$\partial^{\infty}\psi(\bar x):=\{x^*\in X^*\mid (x^*,0)\in N_{{\rm epi}\,\psi}(\bar x,\psi(\bar x))\}$$
is the \textit{singular subdifferential} of $\psi$ at $\bar x$. For a set $\Omega\subset X$ and a point $\bar x\in \Omega$, we have
$$N_{\Omega}(\bar x)=\partial\delta_{\Omega}(\bar x)=\partial^{\infty}\delta_{\Omega}(\bar x),$$
where $\delta_{\Omega}(\bar x)$ is the indicator function of $\Omega$; see \cite[Proposition~1.79]{B-M06}. If $\psi$ depends on two variables $x$ and $y$, and $|\psi(\bar x, \bar y)|<\infty$, then $\partial_x\psi(\bar x,\bar y)$ denotes the Mordukhovich subdifferential of $\psi(.,\bar y)$ at $\bar x$. For any $\bar v\in \partial\psi(\bar x)$,
$$\partial^2\psi(\bar x|\bar v)(u):=D^*(\partial\psi)(\bar x|\bar v)(u)\quad (u\in X^{**}=X)$$ is the \textit{limiting second-order subdifferential} (or the generalized Hessian).

A multifunction $G: Y \rightrightarrows X$ is said to be \textit{locally Lipschitz-like} around $(\bar{y}, \bar{x})\in {\rm gph}\,G$ if there exists a constant $\ell > 0$ and neighborhoods $U$ of $\bar{x}$, $V$ of $\bar{y}$ such that
\begin{equation*}\label{Lipschitz-like property} G(y') \cap U \subset G(y) + \ell \|y' - y\| \bar{B}_{X} \quad \forall y, y' \in V, \end{equation*}
where $\bar{B}_{X}$ denotes the closed unit ball in $X$. When $G$ is locally closed around $(\bar{y}, \bar{x})$, the \textit{Mordukhovich criterion} (see \cite{B-M_TAMS_1993}, \cite[Theorem~9.40]{Rock_Wets_1998}, and \cite[Theorem~4.10]{B-M06}) says that $G$ is locally Lipschitz-like around $(\bar{y}, \bar{x})$ if and only if
\begin{equation*}\label{M_criterion}D^{*}G(\bar{y}, \bar{x})(0) = \{0\}. \end{equation*}

For a multifunction $F: X \times Y \rightrightarrows Z$ and a pair $(\bar{x}, \bar{y}) \in X\times Y$ satisfying $0\in F(\bar{x}, \bar{y})$, we say that the \textit{implicit multifunction} $G: Y \rightrightarrows X$ given by
$G(y) = \{x \in X\mid 0 \in F(x, y)\}$ has the \textit{Robinson stability} at $\omega_{0}:= (\bar{x}, \bar{y}, 0)$ if there exist constants $r > 0$, $\gamma > 0$, and neighborhoods $U$ of $\bar{x}$, $V$ of $\bar{y}$ such that
\begin{equation*}\label{Metric_Regularity}d(x, G(y)) \leq rd(0, F(x, y)) \end{equation*}
for any $(x, y)\in U\times V$ with $d(0, F(x, y)) < \gamma$.
Note that the condition $d(0, F(x, y)) < \gamma$ can be omitted if $F$ is \textit{inner semicontinuous} at $(\bar x, \bar y, 0)$; see~\cite{Huyen_Yen_2016}. Note that, in some cases, the Robinson stability of $G$ at $(\bar x,\bar y,0)$ implies its local Lipschitz-likeness around $(\bar y,\bar x)$; see, e.g., \cite{Gfrerer_Mor_SIOPT2016}. For the generalized linear constraint system studied in \cite{Huyen_Yen_2016}, these properties are equivalent. In the sequel, we will see that the regularity conditions in use guarantee for our stationary point set map to have both properties.

Now, let $f_0$ and $F$ be twice continuously differentiable real-valued functions ($C^2$-functions for brevity) defined on the product $\R^n\times \R^d$ of two Euclidean spaces. For every $w\in\R^d$, we consider the parametric optimization problem
\begin{equation*}\label{Optim_Prob}
(P_w)\quad\ \, {\rm Minimize}\ \, f_0(x,w) \ \; {\rm subject\ to}\ x\in \R^n\ \, {\rm and}\ \, F(x,w)\leq 0.
\end{equation*} The constraint set of $(P_w)$ is $C(w):=\{x\in\R^n\mid F(x,w)\leq 0\}$. The stationary point set of $(P_w)$ is defined by
\begin{equation}\label{KKT_point-set0}
S(w)=\{x\in \R^n\mid 0\in  \nabla_{x}f_0(x,w)+N_{C(w)}(x)\}.
\end{equation}
When $w$ varies on $\R^d$, one has a multifunction  $S:\R^d\rightrightarrows\R^n$ with $S(w)$ being calculated by~\eqref{KKT_point-set0}. Setting $f(x,w)=g(F(x,w))=(g\circ F)(x,w)$, where $g(y)=\delta_{\R_-}(y)$, i.e., $g(y)=0$ for $y\in (-\infty,0]$ and $g(y)=+\infty$ for $y>0$, we can rewrite \eqref{KKT_point-set0} as
\begin{equation}\label{KKT_point-set}
S(w)=\{x\in \R^n\mid 0\in  \nabla_{x}f_0(x,w)+\partial_{x}f(x,w)\}.
\end{equation}

Fix a vector $w=\bar w\in \R^d$ and suppose that $\bar x\in S(\bar w)$. Since $(P_{\bar w})$ has a single smooth inequality constraint, the Mangasarian-Fromovitz Constraint Qualification is fulfilled at $\bar x\in C(\bar w)$ if and only if
\begin{equation}\label{MFCQ_condition}
\text{If}\ F(\bar x,\bar w)=0,\ \text{then}\ \nabla_xF(\bar x,\bar w)\neq 0.\quad\quad\tag{\textbf{MFCQ}}
\end{equation}
In what follows, we assume that \textbf{(MFCQ)} is valid.
To study the stability of the stationary point set map $S$ around the $(\bar w,\bar x)$ in ${\rm gph}\,S$, we compute the Mordukhovich and the Fr\' echet coderivatives of the partial subdifferential map $\partial_{x}f:\R^n\times \R^d\rightrightarrows\R^n$. In general, there is no explicit formula for the coderivatives of such maps. However, the results of \cite{LeMo04} provide us with some tools which allow us to estimate the coderivative value $D^*S(\bar w|\bar x)(x')$ for every $x'\in\R^n$.

The fulfillment of MFCQ at $(\bar x,\bar w)$ implies that $g(x,w)=g(F(x,w))$ is a strongly amenable in $x$ at $\bar x$ with compatible parameterization in $w$ at $\bar w$. Then, by \cite[Theorem~10.49]{Rock_Wets_1998}, for $(x,w)$ near $(\bar x,\bar w)$, we have
\begin{equation}\label{Implication_of_s_amenability1}
\partial f(x,w)=\nabla F(x,w)^*(\partial g(F(x,w)))
\end{equation}
and
\begin{equation}\label{Implication_of_s_amenability2}
\partial_xf(x,w)=\nabla_xF(x,w)^*(\partial g(F(x,w)));
\end{equation}
see \cite[formulas~(14) and (15)]{LeMo04}.

In order to estimate the limiting second-order subdifferential of $f$, we need the following result.

\begin{Lemma}\label{Theorem3.1_LeMo04} {\rm (see \cite[Theorem 3.1]{LeMo04})} 
	Suppose that $\bar v\in \partial f(\bar x,\bar w)$. Then, for
	any $v'\in \R^n\times \R^d$,
	\begin{equation*}
	\begin{array}{rl}
	&\partial^{2}f((\bar x,\bar w)|\bar v)(v')\\
	& \subset\ \displaystyle\bigcup_
	{\begin{subarray}{c} \bar y\in \partial g(F(\bar x,\bar w))\; \text{with}\\ \nabla F(\bar x,\bar w)^*\bar y=\bar v \end{subarray}}
	\bigg(\nabla^{2}(\bar {y}\cdot F)(\bar x,\bar w)v'
	+D^*(\partial g\circ F)(\bar x,\bar w)|\bar y)(\nabla F(\bar x,\bar w)v')\bigg),
	\end{array}
	\end{equation*}
	where the function $\bar {y}\cdot F:\R^{n+d}\to \R$ is defined by $(\bar {y}\cdot F)(x,w):=\bar y F(x,w)$. If, in addition, at every $\bar y\in \partial g(F(\bar x,\bar w))$ with $\nabla F(\bar x,\bar w)^*\bar y=\bar v$, one has the second-order constraint qualification
	\begin{equation}\label{Condi17_LeMor}
	\partial^{2}g(F(\bar x,\bar w)|\bar y)(0)\cap {\rm ker}\nabla F(\bar x,\bar w)^*=\{0\},
	\end{equation} then the estimate above for the second-order subdifferential can be refined by replacing the coderivative of the multifunction $\partial g\circ F$ via the inclusion
	\begin{equation*}\label{Upper_esti_coder_compositefunc}
	D^*(\partial g\circ F)((\bar x,\bar w)|\bar y)(\nabla F(\bar x,\bar w)v')\subset \nabla F(\bar x,\bar w)^* \partial^{2}g(F(\bar x,\bar w)|\bar y)(\nabla F(\bar x,\bar w)v').
	\end{equation*}
\end{Lemma}

In our problem $(P_w)$, condition \eqref{Condi17_LeMor} can be omitted. Indeed, $\bar y\in\partial g(F(\bar x,\bar w))$ if and only if $\bar y \in N_{\R_{-}}(F(\bar x,\bar w))$. Hence, $\bar y\geq 0$. Clearly, $${\rm gph}\,\partial g=(\R_{-}\times \{0\})\cup (\{0\}\times \R_{+}).$$
If $F(\bar x,\bar w)<0$, then $\bar y=0$ and $N_{{\rm gph}\,\partial g}(F(\bar x,\bar w),\bar y)=\{0\}\times \R.$
It follows that
\begin{eqnarray*}\begin{array}{rl} \partial^2g(F(\bar x,\bar w)|\,\bar y)(0)& =D^*(\partial g(F(\bar x,\bar w)|\bar y))(0)\\
		& =\{u'\in \R\mid (u',0)\in N_{{\rm gph}\,\partial g}(F(\bar x,\bar w),\bar y)\}=\{0\}.
	\end{array}
\end{eqnarray*}
So \eqref{Condi17_LeMor} is satisfied. If $F(\bar x,\bar w)=0$, then \textbf{(MFCQ)} implies $\nabla F(\bar x,\bar w)\neq 0$. Hence the linear operator $\nabla F(\bar x,\bar w):\R^n\times\R^n\to\R$ is surjective. Thus ${\rm ker}\,\nabla F(\bar x,\bar w)^*=\{0\}$ by~\cite[Lemma~1.18]{B-M06}, and we see that \eqref{Condi17_LeMor} is fulfilled. Therefore, applied to $(P_w)$, Lemma~\ref{Theorem3.1_LeMo04} can be reformulated as follows: \textit{For any} $\bar v\in\partial f(\bar x,\bar w)$ \textit{and} $v'\in \R^n\times \R^d$,
\begin{equation}\label{Upper_esti_2nd_order_subdifferential_of_f}
\partial^{2}f((\bar x,\bar w)|\bar v)(v')\subset\displaystyle\bigcup_
{\begin{subarray}{c} \bar y\in \partial g(F(\bar x,\bar w))\;\text{with}\\ \nabla F(\bar x,\bar w)^*\bar y=\bar v \end{subarray}}
\Big(\nabla^{2}(\bar {y}\cdot F)(\bar x,\bar w)v'
+\Omega_1(\bar y,v')\Big),
\end{equation} \textit{where}
\begin{equation*}\label{Omega_1}
\Omega_1(\bar y,v'):=\nabla F(\bar x,\bar w)^* \partial^{2}g(F(\bar x,\bar w)|\bar y)(\nabla F(\bar x,\bar w)v').
\end{equation*}

\begin{Remark} {\rm Concerning the paper \cite{MR_2012}, 	
		observe that the set $\partial^{2}f((\bar x,\bar w)|\bar v)(v')$ in  formula \eqref{Upper_esti_2nd_order_subdifferential_of_f} is analogous to the set $\widetilde\varphi^2_x(\bar x,\bar w,\bar y)(u)$ (a value of the extended partial second-order subdifferential) in formula~(3.4) of that work. A careful checking shows that equality (3.4) of \cite{MR_2012} implies the upper estimate \eqref{Upper_esti_2nd_order_subdifferential_of_f}.}
\end{Remark}

In what follows, for any $\bar v=(\bar v_x,\bar v_w)\in\R^n\times \R^d$, we put ${\rm proj}_1\bar v=\bar v_x$. The upper estimation for the coderivative values of the stationary point set map $S$ given by Levy and Mordukhovich \cite{LeMo04} requires the following \textit{regularity condition}: \textit{For any} $v'_1\in\R^n$,
\begin{equation}\label{Condi11_LeMor}
0\in \nabla^{2}f_0(\bar x,\bar w)^*(v'_1,0)+\displaystyle\bigcup_{\begin{subarray}{c} \bar v\in \partial f(\bar x,\bar w)\;\text{with}\\ {\rm proj}_1\bar v=-\nabla_{x}f_0(\bar x,\bar w)\end{subarray}}\partial^2f((\bar x,\bar w)|\bar v)(v'_1,0)\ \; \Longrightarrow\ \; v'_1=0
\end{equation}
(see \cite[formula~(11)]{LeMo04}).
For our problem $(P_w)$, by the assumption \textbf{(MFCQ)} and formula~\eqref{Implication_of_s_amenability1}, we have $\partial f(\bar x,\bar w)=\nabla F(\bar x,\bar w)^*(\partial g(\bar x,\bar w))$. In addition, it is easy to show that, for every $\bar y\in\partial g(\bar x,\bar w)$, ${\rm proj}_1\left(\nabla F(\bar x,\bar w)^*\bar y\right)=\nabla_xF(\bar x,\bar w)^*\bar y$. Hence
\begin{equation}\label{equality_Omega_2}\begin{array}{rl}
& \displaystyle\bigcup_{\begin{subarray}{c} \bar v\in \partial f(\bar x,\bar w)\;\text{with}\\ {\rm proj}_1\bar v=-\nabla_{x}f_0(\bar x,\bar w)\end{subarray}}\partial^2f((\bar x,\bar w)|\bar v)(v'_1,0)\\
&=\displaystyle\bigcup_{\begin{subarray}{c} \bar y\in \partial g(F(\bar x,\bar w))\;\text{with}\\\nabla_{x}F(\bar x,\bar w)^*\bar y=-\nabla_{x}f_0(\bar x,\bar w)\end{subarray}}\partial^2f((\bar x,\bar w)|\nabla F(\bar x,\bar w)^*\bar y)(v'_1,0).
\end{array}\end{equation}
So \eqref{Condi11_LeMor} is equivalent to the following condition:
\begin{equation}\label{C0_condition}
0\in \nabla^{2}f_0(\bar x,\bar w)^*(v'_1,0)+\Omega_2(v'_1)\ \; \Longrightarrow\ \; v'_1=0,\quad\quad\tag{\textbf{C0}}
\end{equation}
where \begin{equation}\label{Omega_2} \Omega_2(v'_1):=\displaystyle\bigcup_{\begin{subarray}{c} \bar y\in \partial g(F(\bar x,\bar w))\;\text{with}\\\nabla_{x}F(\bar x,\bar w)^*\bar y=-\nabla_{x}f_0(\bar x,\bar w)\end{subarray}}\partial^2f((\bar x,\bar w)|\nabla F(\bar x,\bar w)^*\bar y)(v'_1,0).\end{equation}

The next result from \cite{LeMo04} provides us with an upper estimation for the values of the coderivative map  $D^*S(\bar w|\bar x): \R^n\rightrightarrows\R^d$.

\begin{Lemma}[{see \cite[Corollary~3.1]{LeMo04}}]\label{Corollary3.1_LeMo04}
	If the regularity condition \eqref{C0_condition} holds then, for each $x'\in\R^n$, the coderivative value $D^*S(\bar w|\bar x)(x')$ is contained in the set of $w'\in \R^d$ for which there exists a vector $v'_1\in \R^n$ with
	\begin{equation*}
	(-x',w')-\nabla^{2}f_0(\bar x,\bar w)^*(v'_1,0)\in \Omega_2(v'_1).
	\end{equation*}
\end{Lemma}

Although it is rather difficult to compute the set $\Omega_2(v'_1)$, we can still estimate it by using \eqref{Upper_esti_2nd_order_subdifferential_of_f}.

\textit{Upper estimates} for the limiting coderivative values of $S$ can be derived from a result of Levy and Mordukhovich \cite[Theorem~2.1]{LeMo04}. But, a constraint qualification must be imposed to have these estimates (see \cite[p.~1020]{LeeYen11AA} for details). Interestingly, due to a result of Lee and Yen \cite[Theorem~3.4]{LeeYen11AA}, sharp \textit{lower estimates} for the Fr\' echet coderivative values of $S$ can be given without any condition. Put $G(x,w)=\nabla_xf_0(x,w)$ and $M(x,w)=\partial_{x}f(x,w)$. Then,
\begin{equation}\label{Sum_for_S} S(w)=\{x\in\R^n\mid 0\in G(x,w)+M(x,w)\}.\end{equation}
Since $\bar x\in S(\bar w)$, $\bar\tau:=(\bar x,\bar w,-\nabla_xf_0(\bar x,\bar w))$ belongs to ${\rm gph}\,M$. Note that ${\rm gph}\,M$ is locally closed around $\bar\tau$. The following result combines the lower estimates with the upper estimates mentioned above.

\begin{Lemma}\label{Combined_estimates} {\rm (see \cite[Theorem~3.4]{LeeYen11AA})}
	The lower estimates
	$$\widehat\Gamma(x')\subset\widehat D^*S(\bar w|\bar x)(x')\subset D^*S(\bar w\mid\bar x)(x'),$$
	where
	\begin{equation*}\label{Widehat_Gamma(x')}
	\widehat\Gamma(x'):=\displaystyle\bigcup_{v'_1\in \R^n}\left\{w'\in \R^d\mid (-x',w')\in \nabla G(\bar x,\bar w)^*v'_1+\widehat D^*M(\bar\tau)(v'_1)\right\},
	\end{equation*} hold for any $x'\in\R^n$. If the constraint qualification
	\begin{equation}\label{C1_condition}
	0\in \nabla G(\bar x,\bar w)^*v'_1+D^*M(\bar\tau)(v'_1)\ \; \Longrightarrow\ \;v'_1=0
	\quad\quad\tag{\textbf{C1}}
	\end{equation}
	is satisfied, then the upper estimate
	$$D^*S(\bar w|\bar x)(x')\subset\Gamma(x'),$$ where
	\begin{equation*}\label{Gamma(x')}
	\Gamma(x'):=\displaystyle\bigcup_{v'_1\in \R^n}\left\{w'\in \R^d\mid (-x',w')\in \nabla G(\bar x,\bar w)^*v'_1+ D^*M(\bar\tau)(v'_1)\right\},
	\end{equation*} is valid for any $x'\in\R^n$. If, in addition, M is graphically regular at $\bar\tau$, then
	$$\widehat\Gamma(x')=\widehat D^*S(\bar w|\bar x)(x')= D^*S(\bar w|\bar x)(x')=\Gamma(x').$$
\end{Lemma}

From Lemma \ref{Combined_estimates}, for any $x'\in\R^n$, $\widehat\Gamma(x')\subset\widehat D^*S(\bar w|\bar x)(x')$.
This implies that $\widehat\Gamma(0)\subset\widehat D^*S(\bar w|\bar x)(0)\subset D^*S(\bar w|\bar x)(0)$. If we put $\widetilde{M}(x,w)=G(x,w)+M(x,w)$, then by the Fr\'echet coderivative sum rule with equalities \cite[Theorem~1.62]{B-M06},
\begin{equation*}\label{Sum_rule_Fr1}
\widehat D^*\widetilde{M}(\omega_0)(v'_1)=\nabla G(\bar x,\bar w)^*v'_1+\widehat D^*M(\bar\tau)(v'_1)
\end{equation*}
for any $v'_1\in\R^n$, where $\omega_0:=(\bar x,\bar w,0)\in {\rm gph}\,\widetilde{M}$.  Therefore, we can write 	$$\widehat\Gamma(x')=\displaystyle\bigcup_{v'_1\in \R^n}\left\{w'\in \R^d\mid (-x',w')\in\widehat D^*\widetilde{M}(\omega_0)(v'_1)\right\}.$$
Note that $0\in \widehat\Gamma(0)$. According to the Mordukhovich criterion, if $S$ is locally Lipschitz-like around $(\bar w,\bar x)$, then $D^*S(\bar w|\bar x)(0)=\{0\}$ and $\widehat\Gamma(0)=\{0\}$ as a result. In addition, if the constraint qualification \eqref{C1_condition} is fulfilled, then Lemma~\ref{Combined_estimates} yields $D^*S(\bar w|\bar x)(x')\subset\Gamma(x')$ for any $x'\in\R^n$.
In particular, $D^*S(\bar w|\bar x)(0)\subset\Gamma(0)$. Hence, if \eqref{C1_condition} is valid and $\Gamma(0)=\{0\}$, then $$D^*S(\bar w|\bar x)(0)=\{0\}.$$ So, due to the Mordukhovich criterion, $S$ is locally Lipschitz-like around $(\bar w,\bar x)$. This idea has been presented in \cite{LeeYen11AA} and we will follow it throughout this paper.

Put $\mathcal{D}=\{(x,w)\in \R^n\times\R^d\mid F(x,w)\leq 0\}$. If $F(\bar x,\bar w)<0$, then $(\bar x, \bar w)$ is an \textit{interior point} of~$\mathcal{D}$. If  $F(\bar x,\bar w)=0$, then $(\bar x, \bar w)$ is a \textit{boundary point} of $\mathcal{D}$. In the next two sections, we will consider separately these two possibilities of the reference point $(\bar x, \bar w)$. Remind that $\bar w\in \R^d$ and $\bar x\in S(\bar w)$ are fixed and all the notations of this section are kept unchanged.

\section{The Robinson Stability of the Stationary Point Set Map}

Now we turn attention to the Robinson stability of the stationary point set map $S$ of the problem $(P_w)$. As in the preceding sections, we assume the fulfillment of the condition \textbf{(MFCQ)}, which requires that  $\nabla_xF(\bar x,\bar w)\neq 0$ whenever $F(\bar x,\bar w)=0$.

From \cite[Theorem~1.62]{B-M06}, we have a formula similar to \eqref{Sum_rule_Fr1}:
\begin{equation}\label{Limiting_sum_rule}
D^*\widetilde{M}(\omega_0)(v'_1)=\nabla G(\bar x,\bar w)^*v'_1+ D^*M(\bar\tau)(v'_1)
\end{equation}for any $v'_1\in\mathbb{R}^n$. So, condition \eqref{C1_condition} can be rewritten as
$${\rm ker}\,D^*\widetilde{M}(\omega_0)=\{0\}.$$
By \cite[Theorem~3.1]{Yen_Yao_NA_2009}, $S$ has the Robinson stability at $\omega_0=(\bar x,\bar w,0)\in{\rm gph}\,\widetilde{M}$ if \eqref{C1_condition} and the condition
\begin{equation}\label{C2_condition}
\left\{w'\in \mathbb{R}^d\mid\exists v'_1\in \mathbb{R}^n\ {\rm with}\ (0,w')\in D^*\widetilde{M}(\omega_0)(v'_1)\right\}=\{0\},\quad\quad\tag{\textbf{C2}}
\end{equation} is fulfilled. By \eqref{Limiting_sum_rule} we can rewrite \eqref{C2_condition} equivalently as
$$\left\{w'\in \mathbb{R}^d\mid\exists v'_1\in \mathbb{R}^n\ {\rm with}\ (0,w')\in \nabla G(\bar x,\bar w)^*v'_1+ D^*M(\bar\tau)(v'_1)\right\}=\{0\}.$$
With $\Gamma(x')$ defined by  \eqref{Gamma(x')}, we can assert that \eqref{C2_condition} is equivalent to the requirement $\Gamma(0)=\{0\}$. In the proof of \cite[Corollary~2.2]{LeMo04}, the authors have commented that the constraint qualification \eqref{Condi11_LeMor}, which is equivalent to \eqref{C0_condition}, is stronger than \eqref{C1_condition}. Now we go back to three cases considered in Sects.~{\color{blue}3} and~{\color{blue}4} of Part 1.

First, for the case $(\bar x,\bar w)\in {\rm int}\,\mathcal{D}$, we have shown in Sect.~{\color{blue}3} of Part 1 that $D^*M(\bar\tau)(v'_1)=\{0\}$ for any $v'_1\in\mathbb{R}^n$. So, condition \eqref{C2_condition} becomes
$$\left\{w'\in \mathbb{R}^d\mid\exists v'_1\in \mathbb{R}^n\ {\rm with}\ (0,w')\in \nabla G(\bar x,\bar w)^*v'_1\right\}=\{0\}.$$
Since $\nabla G(\bar x,\bar w)^*v'_1=(\nabla^2_{xx} f_0(\bar x,\bar w)v'_1,\nabla^2_{wx} f_0(\bar x,\bar w)v'_1)$, this is equivalent to
$$\left\{w'\in \mathbb{R}^d\mid \exists v'_1\in \mathbb{R}^n\ {\rm with}\ \nabla^{2}_{xx}f_0(\bar x,\bar w)v'_1=0,\; w'=\nabla^{2}_{wx}f_0(\bar x,\bar w)v'_1\right\}=\{0\}.$$
The latter can be rewritten as
	\begin{equation}\label{Condi_tildeC2}
	{\rm ker}\,\nabla^{2}_{xx}f_0(\bar x,\bar w)\subset{\rm ker}\,\nabla^{2}_{wx}f_0(\bar x,\bar w).
	\end{equation}
 Besides, if the condition
 \begin{equation}\label{Condi_tildeC1_S}
 {\rm ker}\,\nabla^{2}_{xx}f_0(\bar x,\bar w)\cap {\rm ker}\,\nabla^{2}_{wx}f_0(\bar x,\bar w)=\{0\}
 \end{equation}
 is fulfilled, then \eqref{C0_condition} is valid and \eqref{C1_condition} is also valid as a result. Thus, \textit{if \eqref{Condi_tildeC1_S} and \eqref{Condi_tildeC2} are simultaneously satisfied, then $S$ has the Robinson stability at $\omega_0$.}

Let us move to the next case where $F(\bar x,\bar w)=0$ and the Lagrange multiplier $\lambda$ corresponding to the stationary point $\bar x\in S(\bar w)$ is positive. First, it is worth to stress that for $(P_w)$, the assumptions (i), (ii), and (10) in \cite[Proposition~2.1]{LeMo04} are fulfilled. So, from \cite[Corollary~2.1]{LeMo04}, for any $v'_1\in\mathbb{R}^n$,
$$D^*M(\bar \omega)(v'_1)\subset \displaystyle\bigcup_{\begin{subarray}{c} \bar v\in \partial f(\bar x,\bar w)\;\text{with}\\ {\rm proj}_1\bar v=-\nabla_{x}f_0(\bar x,\bar w)\end{subarray}}\partial^2f((\bar x,\bar w)|\bar v)(v'_1,0).$$  With $\Omega_2(v'_1)$ defined by \eqref{Omega_2}, using \eqref{equality_Omega_2} we have
$$\Omega_2(v'_1)=\displaystyle\bigcup_{\begin{subarray}{c} \bar v\in \partial f(\bar x,\bar w)\;\text{with}\\ {\rm proj}_1\bar v=-\nabla_{x}f_0(\bar x,\bar w)\end{subarray}}\partial^2f((\bar x,\bar w)|\bar v)(v'_1,0).$$
Hence,
$D^*M(\bar \omega)(v'_1)\subset \Omega_2(v'_1)$  for any $v'_1\in\mathbb{R}^n$. Therefore, from formula  \eqref{Gamma(x')} and the presentation $\nabla G(\bar x,\bar w)^*(v'_1)=\nabla^{2}f_0(\bar x,\bar w)^*(v'_1,0)$, we have
\begin{equation}\label{Upper_esti_of_Gamma(x')}
\Gamma(x')\subset\left\{w'\in\mathbb{R}^d\mid\exists v'_1\in \mathbb{R}^n\ {\rm with}\ (-x',w')-\nabla^{2}f_0(\bar x,\bar w)^*(v'_1,0)\in \Omega_2(v'_1)\right\}
\end{equation} for any $x'\in\mathbb{R}^n$. This implies that  $\Gamma(x')$ is contained in $\Gamma_2(x')$ which is defined in Subsect.~{\color{blue}4.1} of Part 1. In particular, $\Gamma(0)\subset\Gamma_2(0)$. So, if $\Gamma_2(0)=\{0\}$, then  $\Gamma(0)=\{0\}$. We have shown that $\Gamma_2(0)=\{0\}$ if and only if the inclusion
\begin{equation}\label{Condi_tildeC2_haty>0}
{\rm ker}\,A_1\cap \left({\rm ker}\,\nabla_{x}F(\bar x,\bar w)\times \R\right)\subset {\rm ker}\,A_2.
\end{equation}
is valid. Therefore, if \eqref{Condi_tildeC2_haty>0} is satisfied, then $\Gamma(0)=\{0\}$  which implies the fulfillment of \eqref{C2_condition}.
Let 
\begin{equation}\label{matrix_A1} A_1=\left[\begin{array}{ccc} \nabla^{2}_{xx}f_0(\bar x,\bar w)+\lambda\nabla^{2}_{xx}F(\bar x,\bar w) &\ \nabla_{x}F(\bar x,\bar w) \end{array}\right]\in \R^{n\times (n+1)}\end{equation}
and
\begin{equation}\label{matrix_A2} A_2=\left[\begin{array}{ccc} \nabla^{2}_{wx}f_0(\bar x,\bar w)+\lambda\nabla^{2}_{wx}F(\bar x,\bar w) &\ \nabla_{w}F(\bar x,\bar w)\end{array}\right]\in \R^{d\times (n+1)},\end{equation}
where $\nabla_{x}F(\bar x,\bar w)$ and  $\nabla_{w}F(\bar x,\bar w)$ are interpreted as column vectors. If the equality
\begin{equation}\label{Condi_tildeC1_S_haty>0}
{\rm ker}\,A_1\cap {\rm ker}\,A_2\cap \left({\rm ker}\,\nabla_{x}F(\bar x,\bar w)\times \R\right)=\{(0,0)\}.
\end{equation}
is satisfied then, as shown in Subsect.~{\color{blue}4.1} of Part 1, \eqref{C0_condition} is fulfilled; consequently, \eqref{C1_condition} is valid. Thus, \textit{in the case under our consideration, once \eqref{Condi_tildeC1_S_haty>0} and \eqref{Condi_tildeC2_haty>0} are simultaneously satisfied, $S$ has the Robinson stability at $\omega_0$.}

Finally, we consider the case where $F(\bar x,\bar w)=0$ and the Lagrange multiplier $\lambda$ corresponding to the stationary point $\bar x\in S(\bar w)$ equals to zero. In this case, if 
	\begin{equation}\label{Condi_tildeC1_S_haty=0}
	{\rm ker}\,A'_1\cap {\rm ker}\,A'_2=\{0\},
	\end{equation}
	where 
		$$A'_1:=\left[\begin{array}{ccc} \nabla^{2}_{xx}f_0(\bar x,\bar w) &\  \nabla_{x}F(\bar x,\bar w)  \end{array}\right]\in \R^{n\times (n+1)},$$
		and
		$$A'_2:=\left[\begin{array}{ccc} \nabla^{2}_{wx}f_0(\bar x,\bar w) &\ \nabla_{w}F(\bar x,\bar w)  \end{array}\right]\in \R^{d\times (n+1)},$$ 
 is valid, then \eqref{C0_condition} holds (see Subsect.~{\color{blue}4.2} of Part 1). So, \eqref{Condi_tildeC1_S_haty=0} guarantees the validity of \eqref{C1_condition}. Concerning condition \eqref{C2_condition}, we will show that if the conditions 
 \begin{equation}\label{CondiC2a_haty=0}
 {\rm ker}\,A'_1\cap \left({\rm ker}\,\nabla_{x}F(\bar x,\bar w)\times \R\right)\subset {\rm ker}\,A'_2,
 \end{equation}
 \begin{equation}\label{CondiC2b_haty=0}
 {\rm ker}\,A'_1\cap\Delta_1\subset {\rm ker}\,A'_2,
 \end{equation} and
 \begin{equation}\label{CondiC2c_haty=0}
 {\rm ker}\,\nabla^{2}_{xx}f_0(\bar x,\bar w)\cap\Delta_2\subset{\rm ker}\,\nabla^{2}_{wx}f_0(\bar x,\bar w).
 \end{equation}
are satisfied, then \eqref{C2_condition} is fulfilled.  

Let $\Gamma_3(x')$ be the set of vectors $w'\in\R^d$ for which there exists $v'_1\in\R^n$ with
\begin{equation*}
\begin{cases}\left(-x'-\nabla^{2}_{xx}f_0(\bar x,\bar w)v'_1,w'-\nabla^{2}_{wx}f_0(\bar x,\bar w)v'_1\right)\in \{\gamma\nabla F(\bar x,\bar w)\mid \gamma\in\R\},\\
\nabla_{x}F(\bar x,\bar w)v'_1=0,
\end{cases}
\end{equation*} or
\begin{equation*}
\begin{cases}\left(-x'-\nabla^{2}_{xx}f_0(\bar x,\bar w)v'_1,w'-\nabla^{2}_{wx}f_0(\bar x,\bar w)v'_1\right)\in \{\gamma\nabla F(\bar x,\bar w)\mid \gamma\in\R_+\},\\
\nabla_{x}F(\bar x,\bar w)v'_1>0,
\end{cases}
\end{equation*} or
\begin{equation*}
\begin{cases}
-x'-\nabla^{2}_{xx}f_0(\bar x,\bar w)v'_1=0,\ \, w'-\nabla^{2}_{wx}f_0(\bar x,\bar w)v'_1=0,\\
\nabla_{x}F(\bar x,\bar w)v'_1<0.
\end{cases}
\end{equation*}
As it has been proved in Subsect.~{\color{blue}4.2} of Part 1,
\begin{equation}\label{Omega2_degenerate}
\Omega_2(v'_1)\subset
\begin{cases}
\{\gamma\nabla F(\bar x,\bar w)\mid \gamma\in\R\}& \mbox{ if } \nabla_{x}F(\bar x,\bar w)v'_1=0,\\
\{\gamma\nabla F(\bar x,\bar w)\mid \gamma\in\mathbb{R_{+}}\}& \mbox{ if } \nabla_{x}F(\bar x,\bar w)v'_1>0,\\
\{0\} &  \mbox{ if } \nabla_{x}F(\bar x,\bar w)v'_1<0.
\end{cases}
\end{equation}

From  \eqref{Upper_esti_of_Gamma(x')} and the inclusion \eqref{Omega2_degenerate} we have $\Gamma(x')\subset \Gamma_3(x')$ for any $x'\in\mathbb{R}^n$. In particular, $\Gamma(0)\subset \Gamma_3(0)$. Hence, if $\Gamma_3(0)=\{0\}$, then $\Gamma(0)=\{0\}$. We have shown that $\Gamma_3(0)=\{0\}$ holds if and only if the system \eqref{CondiC2a_haty=0}--\eqref{CondiC2c_haty=0} is satisfied. Thus, the validity of \eqref{CondiC2a_haty=0}--\eqref{CondiC2c_haty=0} implies $\Gamma(0)=\{0\}$ which yields the fulfillment of \eqref{C2_condition}. Therefore, \textit{if \eqref{Condi_tildeC1_S_haty=0} and the system \eqref{CondiC2a_haty=0}--\eqref{CondiC2c_haty=0} are simultaneously satisfied, then $S$ has the Robinson stability at $\omega_0$.}

We have thus shown that the sufficient conditions for $S$ being locally Lipschitz-like around $(\bar w,\bar x)$ in each case also guarantee for $S$ having the Robinson stability at $\omega_0$.

Our results on the Robinson stability of $S$ are summarized as follows.

\begin{Theorem}\label{Thm_Robinson_stability} The stationary point set map $S$ of $(P_w)$ has the Robinson stability at $\omega_0=(\bar x,\bar w,0)$ if one of the following is valid:
	\begin{description}
		\item[{\rm (a)}] $F(\bar x,\bar w)<0$ and the condition \begin{equation}\label{Aubin_Suffi_Condi1}
		{\rm ker}\,\nabla^{2}_{xx}f_0(\bar x,\bar w)=\{0\}
		\end{equation} holds;
		\item[{\rm (b)}] $F(\bar x,\bar w)=0$, the Lagrange multiplier $\lambda$ corresponding to the stationary point $\bar x\in S(\bar w)$ is positive, and 
		\begin{equation}\label{Aubin_Suffi_Condi2}
		{\rm ker}\,A_1\cap \left({\rm ker}\,\nabla_{x}F(\bar x,\bar w)\times \R\right)=\{0\};
		\end{equation}
		\item[{\rm (c)}] $F(\bar x,\bar w)=0$, the Lagrange multiplier $\lambda$ corresponding to the stationary point $\bar x\in S(\bar w)$ equals to zero, and \begin{equation}\label{Aubin_Suffi_Condi3}
		\begin{cases}
		{\rm ker}\,A'_1\cap {\rm ker}\,A'_2=\{0\},\\
		{\rm ker}\,A'_1\cap \left({\rm ker}\,\nabla_{x}F(\bar x,\bar w)\times \R\right)\subset {\rm ker}\,A'_2,\\
		{\rm ker}\,A'_1\cap\Delta_1\subset {\rm ker}\,A'_2,\\
		{\rm ker}\,\nabla^{2}_{xx}f_0(\bar x,\bar w)\cap\Delta_2\subset{\rm ker}\,\nabla^{2}_{wx}f_0(\bar x,\bar w).
		\end{cases}
		\end{equation} 
	\end{description}
\end{Theorem}

It is worthy to stress that the Robinson stability of $S$ at $\omega_0$ is available for the examples of the previous section where our sufficient conditions for the local Lipschitz-likeness of $S$ around $(\bar w,\bar x)$ are fulfilled.

\section{Applications to Quadratic Programming}

In this section, the above general results are applied to a class of nonconvex quadratic programming problems. Namely, we will consider the problems of minimizing a linear-quadratic function under one linear-quadratic functional constraint. Special cases of such problems have been considered, e.g., in \cite{LeeTamYen_SIOPT2012}, \cite{LeeYen14NA}, and \cite{QuiYen_SIOPT2014}.

Denote by $\mathcal{S}_n$ the space of $n\times n$ symmetric matrices. Let $D, A\in \mathcal{S}_n$, $c$ and $b$ be vectors in $\mathbb{R}^n$, and $\alpha$ a real number. Put $w=(w_1,w_2)$ with $w_1:=(D,c)$ and $w_2:=(A,b,\alpha)$. Denote the problem $(P_w)$ with
$f_0(x,w)=\frac{1}{2}x^TDx+c^Tx$ and $F(x,w)=\frac{1}{2}x^TAx+b^Tx+\alpha$ by $(QP_w)$. For convenience, we put $W_1=\mathcal{S}_n \times\mathbb{R}^n$, $W_2=\mathcal{S}_n\times\mathbb{R}^n\times\mathbb{R}$, and $W=W_1\times W_2$.
Fix a vector $\bar w=(\bar w_1,\bar w_2)\in W$ with $\bar w_1=(\bar D,\bar c)$, $\bar w_2=(\bar A,\bar b,\bar\alpha)$, and suppose that a stationary point $\bar x\in S(\bar w)$ is given.

To ease the description of certain second order differential operators, sometimes we will present the matrices $D$ and $A$ in the following column forms
\begin{equation*}\label{Column_forms}
D=\left(\begin{matrix}
d^T_{1}\\d^T_2\\\vdots\\d^T_n
\end{matrix}\right),\quad A=\left(\begin{matrix}
a^T_{1}\\a^T_2\\\vdots\\a^T_n
\end{matrix}\right),
\end{equation*}
where $d_i=(d_{i1}\dots d_{in})$ and $a_i=(a_{i1}\dots a_{in})$ are, respectively, the $i$-th row of $D$ and the $i$-th row of $A$. We have $\nabla_xf_0(\bar x,\bar w)=\bar D\bar x+\bar c$,
$$\nabla_{w_1}f_0(\bar x,\bar w)=
\left(\begin{matrix}
\frac{1}{2}\bar x_1\bar x_1\ \dots\ \frac{1}{2}\bar x_1\bar x_n\ \dots\ \frac{1}{2}\bar x_n\bar x_1\ \dots\ \frac{1}{2}\bar x_n\bar x_n\ \bar x_1\ \dots\ \bar x_n
\end{matrix}\right)^T,$$
$\nabla^2_{xx}f_0(\bar x,\bar w)=\bar D$, $\nabla^2_{w_2x}f_0(\bar x,\bar w)=0_{W_2}$, and
$$\nabla^2_{w_1x}f_0(\bar x,\bar w)=
\begin{pmatrix}
\bar X & \dots & 0\\
&\ddots\\
0 & \dots & \bar X\\
1 & \dots & 0\\
&\ddots\\
0 & \dots & 1
\end{pmatrix}.$$
Here, $\bar X:=\begin{pmatrix} \bar x_1\\ \vdots\\ \bar x_n \end{pmatrix}$ is an $n\times 1$ matrix.
Similarly, $\nabla_xF(\bar x,\bar w)=\bar A\bar x+\bar b$, $\nabla^2_{xx}F(\bar x,\bar w)=\bar A$,
$$\nabla_{w_2}F(\bar x,\bar w)=
\left(\begin{matrix}
\frac{1}{2}\bar x_1\bar x_1\ \dots\ \frac{1}{2}\bar x_1\bar x_n\ \dots\ \frac{1}{2}\bar x_n\bar x_1\ \dots\ \frac{1}{2}\bar x_n\bar x_n\ \bar x_1\ \dots\ \bar x_n\ 1
\end{matrix}\right)^T,$$
$\nabla^2_{w_1x}F(\bar x,\bar w)=0_{W_1}$, and
$$\nabla^2_{w_2x}F(\bar x,\bar w)=
\begin{pmatrix}
\bar X & \dots & 0\\
&\ddots\\
0 & \dots & \bar X\\
1 & \dots & 0\\
&\ddots\\
0 & \dots & 1\\
0 & \dots & 0
\end{pmatrix}.$$
We have
$$\nabla^2_{wx}f_0(\bar x,\bar w)=\left(\begin{matrix}
\nabla^2_{w_1x}f_0(\bar x,\bar w)\ \ \nabla^2_{w_2x}f_0(\bar x,\bar w)
\end{matrix}\right).$$
Since $\nabla^2_{w_2x}f_0(\bar x,\bar w)=0$,
$$ {\rm ker}\,\nabla^2_{wx}f_0(\bar x,\bar w)=\left\{v'_1\in\mathbb{R}^n\mid \nabla^2_{w_1x}f_0(\bar x,\bar w)v'_1=0\right\}=\{0\}.$$

First, we consider the case of interior points $(\bar x,\bar w)$, i.e., $F(\bar x,\bar w)<0$. The conditions \eqref{Condi_tildeC1_S}, \eqref{Condi_tildeC2}, and \eqref{Aubin_Suffi_Condi1} are equivalent due to ${\rm ker}\,\nabla^2_{wx}f_0(\bar x,\bar w)=\{0\}$. Thus, by Theorem~3.1 of Part 1, the stationary point set map $S$ of $(P_w)$ is locally Lipschitz-like around $(\bar w,\bar x)$ if and only if ${\rm ker}\,\nabla^2_{xx}f_0(\bar x,\bar w)=\{0\}$, or ${\rm ker}\,\bar D=\{0\}$. In other words, $S$ is locally Lipschitz-like around $(\bar w,\bar x)$ if and only if matrix $\bar D$ is nonsingular. By Theorem~\ref{Thm_Robinson_stability}, this condition is sufficient for $S$ having the Robinson stability at $\omega_0$.

Next, consider the second case where $(\bar x,\bar w)$ is a boundary point of $\mathcal{D}$ and the Lagrange multiplier~$\lambda$ corresponding to $\bar x\in S(\bar w)$ is positive. As in Part~1, $\lambda$ is defined by \begin{equation*}\label{Lagrange_multiplier}
\nabla_{x}f_0(\bar x,\bar w)+\lambda\nabla_{x}F(\bar x,\bar w)=0,
\end{equation*}  which is rewritten as
\begin{equation}\label{QP_Lagrange_multiplier}
\lambda(\bar A\bar x+\bar b)=-(\bar D\bar x+\bar c).
\end{equation} We have
$\nabla^{2}_{xx}f_0(\bar x,\bar w)+\lambda\nabla^{2}_{xx}F(\bar x,\bar w)=\bar D+\lambda\bar A$	 and
$$\nabla^{2}_{wx}f_0(\bar x,\bar w)+\lambda\nabla^{2}_{wx}F(\bar x,\bar w)\\
=\begin{pmatrix}
\bar X & \dots & 0\\
&\ddots\\
0 & \dots & \bar X\\
1 & \dots & 0\\
&\ddots\\
0 & \dots & 1\\
\lambda\bar X & \dots & 0\\
&\ddots\\
0 & \dots & \lambda\bar X\\
\lambda & \dots & 0\\
&\ddots\\
0 & \dots & \lambda\\
0 & \dots & 0
\end{pmatrix}.$$
Now, the matrices $A_1$ and $A_2$ defined in Sect.~{\color{blue}3} are described as follows
$$A_1=\left[\begin{array}{ccc} \bar D+\lambda\bar A & \bar A\bar x+\bar b\end{array}\right]$$
and
$$A_2=\begin{pmatrix}
\bar X & \dots & 0 & 0\\
&\ddots&&\vdots\\
0 & \dots & \bar X & 0\\
1 & \dots & 0 & 0\\
&\ddots&&\vdots\\
0 & \dots & 1 & 0\\
\lambda\bar X & \dots & 0 & \frac{1}{2}\bar x_1\bar x_1\\
&\ddots&&\vdots\\
0 & \dots & \lambda\bar X & \frac{1}{2}\bar x_1\bar x_n\\
\lambda & \dots & 0 & \bar x_1\\
&\ddots&&\vdots\\
0 & \dots & \lambda & \bar x_n\\
0 & \dots & 0 & 1
\end{pmatrix}.$$
Hence, ${\rm ker}\, A_2=\{0\}$. This implies that
\eqref{Condi_tildeC1_S_haty>0} is automatically satisfied. So, according to Theorem~4.1 of Part 1, $S$ is locally Lipschitz-like around $(\bar w,\bar x)$ if and only if \eqref{Condi_tildeC2_haty>0} is fulfilled. Note that
$${\rm ker}\,A_1=\{(v'_1,\gamma)\in\mathbb{R}^n\times\mathbb{R}\mid (\bar D+\lambda\bar A)v'_1+\gamma(\bar A\bar x+\bar b)=0\}$$
and
$${\rm ker\,}\nabla_xF(\bar x,\bar w)=\{v'_1\in\mathbb{R}^n\mid (\bar A\bar x+\bar b)^Tv'_1=0\}.$$
Hence, \eqref{Condi_tildeC2_haty>0} holds if and only if
\begin{equation*}
\begin{cases}
(\bar D+\lambda\bar A)v'_1+\gamma(\bar A\bar x+\bar b)=0\\
(\bar A\bar x+\bar b)^Tv'_1=0\\
\end{cases}
\qquad\ \Longrightarrow\ \,
\begin{cases}
v'_1=0\\\gamma=0
\end{cases}
\end{equation*} or, equivalently,
\begin{equation}\label{QP_haty>0}
\det\begin{pmatrix} \bar D+\lambda\bar A & \bar A\bar x+\bar b\\
(\bar A\bar x+\bar b)^T & 0
\end{pmatrix}\neq 0.
\end{equation}
Thus, $S$ is locally Lipschitz-like around $(\bar w,\bar x)$ if and only if \eqref{QP_haty>0} is satisfied. Moreover, by Theorem~\ref{Thm_Robinson_stability}, \eqref{QP_haty>0} guarantees for $S$ having the Robinson stability at $\omega_0$.

Let us consider the last case where $(\bar x,\bar w)$ is a boundary point of ${\mathcal D}$ and the Lagrange multiplier~$\lambda$ corresponding to $\bar x\in S(\bar w)$ equals to zero. The matrices $A'_1$ and $A'_2$ defined in  Sect.~{\color{blue}3} are described as $A'_1=\left[\begin{array}{ccc} \bar D & \bar A\bar x+\bar b\end{array}\right]$
and
$$A'_2=\begin{pmatrix}
\bar X & \dots & 0 & 0\\
&\ddots&&\vdots\\
0 & \dots & \bar X & 0\\
1 & \dots & 0 & 0\\
&\ddots&&\vdots\\
0 & \dots & 1 & 0\\
0 & \dots & 0 & \frac{1}{2}\bar x_1\bar x_1\\
&\ddots&&\vdots\\
0 & \dots & 0 & \frac{1}{2}\bar x_1\bar x_n\\
0 & \dots & 0 & \bar x_1\\
&\ddots&&\vdots\\
0 & \dots & 0 & \bar x_n\\
0 & \dots & 0 & 1
\end{pmatrix}.$$
Since ${\rm ker}\,A'_2=\{0\}$, using the equality  ${\rm ker}\,\nabla^2_{wx}f_0(\bar x,\bar w)=\{0\}$ we can rewrite \eqref{Aubin_Suffi_Condi3} as
$$\begin{cases}
{\rm ker}\,A'_1\cap \left({\rm ker}\,\nabla_{x}F(\bar x,\bar w)\times \mathbb{R}\right)=\{0\},\\
{\rm ker}\,A'_1\cap\Delta_1=\{0\},\\
{\rm ker}\,\nabla^{2}_{xx}f_0(\bar x,\bar w)\cap\Delta_2=\{0\}.
\end{cases}$$
This condition holds if and only if the following conditions are simultaneously satisfied:
\begin{equation*}
\begin{cases}
\bar Dv'_1+\gamma(\bar A\bar x+\bar b)=0\\
(\bar A\bar x+\bar b)^Tv'_1=0\\
\end{cases}\qquad\ \Longrightarrow\ \,
\begin{cases}
v'_1=0\\\gamma=0,
\end{cases}
\end{equation*}
\begin{equation*}
\begin{cases}
\bar Dv'_1+\gamma(\bar A\bar x+\bar b)=0\\
(\bar A\bar x+\bar b)^Tv'_1>0,\ \gamma\geq 0
\end{cases}
\qquad\ \Longrightarrow\ \,
\begin{cases}
v'_1=0\\\gamma=0,
\end{cases}
\end{equation*}
and
\begin{equation*}
\begin{cases}
\bar Dv'_1=0\\
(\bar A\bar x+\bar b)^Tv'_1<0
\end{cases}
\qquad\ \Longrightarrow\ \, v'_1=0.
\end{equation*} These implications can be rewritten respectively as
\begin{equation}\label{QP_sufficient_condi_haty=0_a}
\det\begin{pmatrix} \bar D & \bar A\bar x+\bar b\\
(\bar A\bar x+\bar b)^T & 0
\end{pmatrix}\neq 0,
\end{equation}
\begin{equation}\label{QP_sufficient_condi_haty=0_b}
[\bar Dv'_1+\gamma(\bar A\bar x+\bar b)=0,\
\gamma\geq 0]
\ \; \Longrightarrow\ \; (\bar A\bar x+\bar b)^Tv'_1\leq 0,
\end{equation}
and
\begin{equation}\label{QP_sufficient_condi_haty=0_c}
\bar Dv'_1=0\ \; \Longrightarrow\ \;
(\bar A\bar x+\bar b)^Tv'_1=0.
\end{equation}
Thus, in accordance with Theorem~4.2 of Part 1, $S$ is locally Lipschitz-like around $(\bar w,\bar x)$ if  \eqref{QP_sufficient_condi_haty=0_a}--\eqref{QP_sufficient_condi_haty=0_c} are satisfied. Moreover, by Theorem~\ref{Thm_Robinson_stability}, the fulfillment of \eqref{QP_sufficient_condi_haty=0_a}--\eqref{QP_sufficient_condi_haty=0_c} is sufficient for $S$ having the Robinson stability at $\omega_0$. Let us consider the necessary condition 
\begin{equation*}\label{Necessary_condi_degenerate}
{\rm ker}\,A'_1\cap\Delta_3\subset {\rm ker}\,A'_2.
\end{equation*}
for the local Lipschitz-like property of $S$  around $(\bar w,\bar x)$, which is now reduced to ${\rm ker}\,A'_1\cap\Delta_3=\{0\}.$
Clearly, this condition is equivalent to
\begin{equation}\label{QP_necessary_condi_haty=0}
\begin{cases}
\bar Dv'_1+\gamma(\bar A\bar x+\bar b)=0\\
(\bar A\bar x+\bar b)^Tv'_1\geq 0,\ \gamma\geq 0\\
\end{cases}\qquad\ \Longrightarrow\ \,
\begin{cases}
v'_1=0\\\gamma=0.
\end{cases}
\end{equation}
By \cite[Theorem~4.1]{Huyen_Yao_Yen_2018a}, \eqref{QP_necessary_condi_haty=0} is a necessary condition for $S$ being locally Lipschitz-like around $(\bar w,\bar x)$.

The obtained results can be formulated as follows.

\begin{Theorem}\label{Thm_QP} The following assertions are true:
	\item[{\rm (a)}] If $F(\bar x,\bar w)<0$, then $S$ is locally around $(\bar w,\bar x)$ if and only if $\det\bar D\neq 0$. Moreover, under this condition, $S$ has the Robinson stability at $\omega_0$;
	\item[{\rm (b)}] If $F(\bar x,\bar w)=0$ and the Lagrange multiplier~$\lambda$ corresponding to $\bar x\in S(\bar w)$ is positive, then $S$ is locally Lipschitz-like around $(\bar w,\bar x)$ if and only if \eqref{QP_haty>0} is fulfilled. This condition is sufficient for $S$ having the Robinson stability at $\omega_0$;
	\item[{\rm (c)}] If $F(\bar x,\bar w)=0$ and the Lagrange multiplier~$\lambda$ corresponding to $\bar x\in S(\bar w)$ is zero, then \eqref{QP_necessary_condi_haty=0} is necessary for $S$ being locally Lipschitz-like around $(\bar w,\bar x)$. Meanwhile, the fulfillment of \eqref{QP_sufficient_condi_haty=0_a}--\eqref{QP_sufficient_condi_haty=0_c} is sufficient for the local Lipschitz-like property of $S$ around $(\bar w,\bar x)$, as well as for the Robinson stability of $S$ at $\omega_0$.
\end{Theorem}

To show how these results can work, we revisit some examples from \cite{QuiYen_SIOPT2014}.

\begin{Example} {\rm (see \cite[Example~4.1]{QuiYen_SIOPT2014})}\label{QP_Exam1}
	{\rm Consider the problem $(QP_w)$ in the case $n=2$. Choosing $\bar w=(\bar D,\bar c,\bar A,\bar b,\bar \alpha)$ with
	\begin{eqnarray*}
		\bar D=\left(\begin{matrix}
			0 &  0\cr
			0 & -8\cr
		\end{matrix}\right),\quad
		\bar c=\left(\begin{matrix}
			1\cr
			0\cr
		\end{matrix}\right)
	\end{eqnarray*}
	and
	\begin{eqnarray*}
		\bar A=\left(\begin{matrix}
			1 &  0\cr
			0 &  1\cr
		\end{matrix}\right),\quad
		\bar b=\left(\begin{matrix}
			0\cr
			0\cr
		\end{matrix}\right),\quad \bar \alpha=-1,
	\end{eqnarray*} one has $f_0(x,\bar w)=-4x_2^2+x_1$ and $F(x,\bar w)=x_1^2+x_2^2-1$.
	To show that $\bar x:=(-\frac{1}{8},\frac{\sqrt{63}}{8})^T$ is a stationary point of $(P_{\bar w})$, we note by \eqref{KKT_point-set} that
	$$S(\bar w)=\left\{x\in\mathbb{R}\mid 0\in \nabla_xf_0(\bar x,\bar w)+\partial_xf(\bar x,\bar w)\right\},$$
	with $f(x,w)=(g\circ F)(x,w)$ and $g(y)=\delta_{\mathbb R_{-}}(y)$ for any $y\in\mathbb{R}$. As $F(\bar x,\bar w)=0$ and $\nabla_xF(\bar x,\bar w)\neq 0$, condition \textbf{(MFCQ)} is valid. So, from \eqref{Implication_of_s_amenability2} we have
	\begin{equation*}
	\begin{array}{rl }
	\partial_xf(\bar x,\bar w)& = \nabla_xF(\bar x,\bar w)^*N_{\mathbb{R}_{-}}(F(\bar x,\bar w))\\
	& =\nabla_xF(\bar x,\bar w)^*\mathbb{R_{+}}\\
	& =\left\{(-\frac{1}{4}\gamma,\frac{\sqrt{63}}{4}\gamma)\mid \gamma\in\mathbb{R}_{+}\right\}.
	\end{array}
	\end{equation*}
	Besides, $\nabla_xf_0(\bar x,\bar w)=(1,-\sqrt{63})^T$. Now, it is clear that $\bar x\in S(\bar w)$. From \eqref{QP_Lagrange_multiplier}, the Lagrange multiplier corresponding to $\bar x$ is $\lambda=8$. Hence,
	\begin{eqnarray*}
		\det\begin{pmatrix} \bar D+\lambda\bar A & \bar A\bar x+\bar b\\
			(\bar A\bar x+\bar b)^T & 0
		\end{pmatrix}=
		\det\left(\begin{matrix}
			8 & 0 &  \frac{1}{8}\cr
			0 & 0 & -\frac{\sqrt{63}}{8}\cr
			-\frac{1}{8} & \frac{\sqrt{63}}{8} &  0\cr
		\end{matrix}\right)
		=\frac{63}{8}.
	\end{eqnarray*}
	So, \eqref{QP_haty>0} is fulfilled. Thus, by Theorem~\ref{Thm_QP}, the stationary point set map $S$ of $(P_w)$ not only is locally Lipschitz-like around $(\bar w,\bar x)$ but also has the Robinson stability at $\omega_0=(\bar x,\bar w,0)$. Similarly, we can show that $\bar x=(-\frac{1}{8},-\frac{\sqrt{63}}{8})^T$ and $\bar x= (-1,0)^T$ belong to $S(\bar w)$ and \eqref{QP_haty>0} is also valid for them.}
\end{Example}

\begin{Example} {\rm (see \cite[Example~4.2]{QuiYen_SIOPT2014})}\label{QP_Exam2}
	{\rm Consider the problem $(QP_w)$ in the case $n=3$ and choose $\bar w=(\bar D,\bar c,\bar A,\bar b,\bar \alpha)$ with
	\begin{eqnarray*}
		\bar D=\left(\begin{matrix}
			0 &  0 & 0\cr
			0 & -8 & 0\cr
			0 & 0 & -8\cr
		\end{matrix}\right),\quad
		\bar c=\left(\begin{matrix}
			1\cr
			0\cr
			0\cr
		\end{matrix}\right)
	\end{eqnarray*}
	and
	\begin{eqnarray*}
		\bar A=\left(\begin{matrix}
			1 &  0 & 0\cr
			0 &  1 & 0\cr
			0 &  0 & 1\cr
		\end{matrix}\right),\quad
		\bar b=\left(\begin{matrix}
			0\cr
			0\cr
			0\cr
		\end{matrix}\right),\quad \bar \alpha=-1.
	\end{eqnarray*} Here, $f_0(x,\bar w)=-4(x_2^2+x_3^2)+x_1$ and $F(x,\bar w)=x_1^2+x_2^2+x_3^2-1$. Arguments similar to those in the previous example show that $\bar x:=(-1,0,0)^T$ is a stationary point of $(P_{\bar w})$ with the associated Lagrange multiplier $\lambda=1$. It is easy to check that \eqref{QP_haty>0} is satisfied. So, by Theorem~\ref{Thm_QP}, the stationary point set map $S$ of $(P_w)$ is  locally Lipschitz-like around $(\bar w,\bar x)$ and it has the Robinson stability at $\omega_0=(\bar x,\bar w,0)$. However, for the stationary points $$\bar x_t:=\left(-\frac{1}{8},\left(\frac{\sqrt{63}}{8}\right)\sin t,\left(\frac{\sqrt{63}}{8}\right)\cos t\right)^T,$$
	with $t\in [0,2\pi)$, which share the common associated Lagrange multiplier $\lambda=8$, \eqref{QP_haty>0} is violated because
	\begin{eqnarray*}
		\det\begin{pmatrix} \bar D+\lambda\bar A & \bar A\bar x+\bar b\\
			(\bar A\bar x+\bar b)^T & 0
		\end{pmatrix}=
		\det\left(\begin{matrix}
			8 & 0 & 0 & \frac{1}{8}\cr
			0 & 0 & 0 & -\frac{\sqrt{63}}{8}\sin t\cr
			0 & 0 & 0 & -\frac{\sqrt{63}}{8}\cos t\cr
			-\frac{1}{8} & \frac{\sqrt{63}}{8}\sin t & \frac{\sqrt{63}}{8}\cos t & 0\cr
		\end{matrix}
		\right)
		=0.
	\end{eqnarray*}
	Thus, by Theorem~\ref{Thm_QP}, $S$ is not locally Lipschitz-like around $(\bar w,\bar x)$.}
\end{Example}	 	

The \textit{parametric trust-region subproblem} (TRS) considered in \cite{LeeTamYen_SIOPT2012,LeeYen11AA,QuiYen_SIOPT2014} is a special case of our quadratic programming problem $(QP_w)$, where $A$ is the unit matrix, $b=0$, and $\alpha<0$.

For (TRS), in the case where $F(\bar x,\bar w)=0$ and the Lagrange multiplier~$\lambda$ corresponding to $\bar x\in S(\bar w)$ is positive, the matrix in \eqref{QP_haty>0} coincides with the matrix $Q(.)$ in \cite[Theorem~5.1]{LeeYen14NA} and \cite[Theorem~4.2]{QuiYen_SIOPT2014}, called the \textit{stability matrix} (see \cite[p.~200]{LeeYen14NA}). Therefore, Theorem~4.2 in~\cite{QuiYen_SIOPT2014}, which only discusses the local Lipschitz-like property, is a consequence of the assertions (a) and (b) of Theorem~\ref{Thm_QP}.

In the case where $F(\bar x,\bar w)=0$ and the Lagrange multiplier~$\lambda$ corresponding to $\bar x\in S(\bar w)$ equals to zero, the matrix in \eqref{QP_sufficient_condi_haty=0_a} coincides with the \textit{stability matrix} $Q_1(.)$ in \cite[Theorem~4.3]{QuiYen_SIOPT2014}. So, condition (4.10) in \cite{QuiYen_SIOPT2014} coincides with our condition \eqref{QP_sufficient_condi_haty=0_a}. The \textit{sufficient condition} for  the local Lipschitz-like property  in \cite[Theorem~4.3]{QuiYen_SIOPT2014} also requires $\det\bar D\neq 0$. Under this assumption, \eqref{QP_sufficient_condi_haty=0_b} and \eqref{QP_sufficient_condi_haty=0_c} are valid if and only if $\bar x^T\bar D^{-1}\bar x\geq 0$. However, our conditions  \eqref{QP_sufficient_condi_haty=0_a}--\eqref{QP_sufficient_condi_haty=0_c} do not require $\det\bar D\neq 0$. Thus, for (TRS), the sufficient conditions in Theorem~\ref{Thm_QP}(c) and in \cite[Theorem~4.3(ii)]{QuiYen_SIOPT2014} are independent results. Finally, note that the \textit{necessary condition} (4.9) in \cite{QuiYen_SIOPT2014} for the local Lipschitz-like property coincides with our condition \eqref{QP_necessary_condi_haty=0}.

\section{Results Obtained by Another Approach}

Following the detailed hints of one referee of this paper, we will compare our results with those which can be obtained by using the theory of strongly regular generalized equations of Robinson \cite{Robinson_1980}.

Suppose that $\bar x\in S(\bar w)$ and the condition {\bf (MFCQ)} is satisfied. It is not difficult to show that, thanks to {\bf (MFCQ)}, there exist a neighborhood $W_0$ of $\bar w$ and a neighborhood $U_0$ of $\bar x$ such that for every $(x,w)\in U_0\times W_0$ one has $N_{C(w)}(x)=\{\lambda\nabla_xF(x,w)\mid \lambda\geq 0\}$ when $F(x,w)=0$ and  $N_{C(w)}(x)=\{0\}$ when $F(x,w)<0$. Hence, for every $(x,w)\in U_0\times W_0$, the condition
$$0\in \nabla_{x}f_0(x,w)+N_{C(w)}(x)$$ is equivalent to the existence of a Lagrange multiplier $\alpha\in \R$ such that
$$0\in 
\left(\begin{matrix}
	\nabla_{x}{\mathcal{L}}(x,\alpha,w)\\
	-F(x,w)
	\end{matrix}\right)
+N_{\R^n\times\R_{+}}(x,\alpha),$$
where ${\mathcal{L}}(x,\alpha,w) := f_0(x,w)+\alpha F(x,w)$. Setting $g(x,\alpha,w)=\left(\begin{matrix}
\nabla_{x}{\mathcal{L}}(x,\alpha,w)\\
-F(x,w)
\end{matrix}\right),$ we consider the \textit{parametric generalized equation} 
\begin{equation}\label{par_gen_eq}
0\in g(x,\alpha,w)+N_{\R^n\times\R_{+}}(x,\alpha)\quad (w\in\R^d)
\end{equation} and denote its solution set by $\widehat S(w)$. Then, $$\widehat S(w)=\{(x,\alpha)\in\R^n\times\R^d\mid 0\in g(x,\alpha,w)+N_{\R^n\times\R_{+}}(x,\alpha)\}$$ and $\widehat S(.)$ is the \textit{implicit multifunction} defined by \eqref{par_gen_eq}. (The writing of the necessary optimality condition of a constrained smooth mathematical programming problem in a form similar to \eqref{par_gen_eq} has been used by Robinson \cite[p.~54]{Robinson_1980}.) From what has been said we have 
\begin{equation}\label{S_via_hat_S} S(w)\cap U_0=\{x\in U_0\mid \exists \alpha\ {\rm s.t.}\ (x,\alpha)\in \widehat S(w)\}\quad (\forall w\in W_0).\end{equation} 

As in Part 1 of this paper and in the preceding sections, we will denote by $\lambda$ the unique multiplier corresponding to $\bar x\in S(\bar w)$. Consider the following three cases.

\medskip
{\sc Case 1}: $F(\bar x,\bar w)<0$. This case has been analyzed in Remark 3.2 of Part 1 of this paper.

\medskip
{\sc Case 2}: $F(\bar x,\bar w)=0$ \textit{and} $\lambda >0$. In accordance with \cite[p.~45]{Robinson_1980}, the unperturbed generalized equation
\begin{equation}\label{par_gen_eq}
0\in g(x,\alpha,\bar w)+N_{\R^n\times\R_{+}}(x,\alpha)
\end{equation} is said to be \textit{strongly regular} at $(\bar x,\lambda)$  if there
exist a constant $\ell_0>0$ and neighborhoods $U$ of the origin in $\R^n\times\R$ and $V$ of $(\bar x,\lambda)$ such that  for every $(x',\alpha')\in U$ one can find a unique vector $(x,\alpha)$ in $V$, denoted by $s_0(x',\alpha')$, satisfying 
\begin{equation*}\label{lin_gen_eq}
\left(\begin{matrix} x'\\ \alpha'\end{matrix}\right)
\in g(\bar x,\lambda,\bar w)+\nabla_{(x,\alpha)}g(\bar x,\lambda,\bar w)((x,\alpha)-(\bar x,\lambda))+N_{\R^n\times\R_{+}}(x,\alpha)
\end{equation*}
and the mapping $s_0:U\to V$ is Lipschitzian on $U$ with
modulus $\ell_0$. Using the condition $\lambda >0$ and the results of Dontchev and Rockafellar \cite{DR_1996}, one can prove next lemma; see Sect.~\ref{appendices} for details. 
\begin{Lemma}\label{new_lemma1}
The generalized equation \eqref{par_gen_eq} is strongly regular at $(\bar x,\lambda)$ iff the matrix
\begin{equation}\label{combined_matrix}\begin{pmatrix} \nabla_{xx}^2{\mathcal{L}}(\bar x,\lambda,\bar w) & \ \; \nabla_{x}F(\bar x,\bar w)\\
\nabla_{x}F(\bar x,\bar w)^T & \ \; 0
\end{pmatrix}\end{equation}
is nonsingular. 
\end{Lemma}

The condition formulated in Lemma \ref{new_lemma1} is equivalent to condition (23) in Part 1, which was renumbered as condition \eqref{Aubin_Suffi_Condi2} in Sect.~3. Indeed, by \eqref{matrix_A1} one has 
\begin{equation*} A_1=\left[\begin{array}{ccc}\nabla_{xx}^2{\mathcal{L}}(\bar x,\lambda,\bar w)  & \; \nabla_{x}F(\bar x,\bar w) \end{array}\right].\end{equation*} Hence, $(x',\tau')\in {\rm ker}\, A_1\cap \big({\rm ker}\,\nabla_{x}F(\bar x,\bar w)\times\R\big)$ iff
$$\begin{pmatrix} \nabla_{xx}^2{\mathcal{L}}(\bar x,\lambda,\bar w) & \ \; \nabla_{x}F(\bar x,\bar w)\\
\nabla_{x}F(\bar x,\bar w)^T & \ \; 0
\end{pmatrix}\begin{pmatrix}x' \\ \tau'\end{pmatrix}=\begin{pmatrix} 0 \\ 0\end{pmatrix}.$$ Thus, the matrix in \eqref{combined_matrix} is nonsingular if \eqref{Aubin_Suffi_Condi2} is valid. Now, applying Theorem~2.1 from \cite{Robinson_1980} to the parametric generalized equation \eqref{par_gen_eq}, we can assert that if \eqref{par_gen_eq} is strongly regular at $(\bar x,\lambda)$, then the implicit multifunction $\widehat S(.)$ has a \textit{single-valued localization} \cite[p.~4]{DR_2009} around $\bar w$ for $(\bar x,\lambda)$ which is  Lipschitz continuous in a neighborhood of $\bar w$. This means that there exist $\ell>0$, a neighborhood $W$ of $\bar w$, a neighborhood $U$ of $\bar x$, and neighborhood $V$ of $\bar x$ such that for each $w\in W$ there is a unique vector $(x(w),\alpha(w))$, denoted by $\hat s(w)$, in $U\times V$ satisfying the equation \eqref{par_gen_eq} and $\|\hat s(w_2)-\hat s(w_1)\|\leq\ell\|w_2-w_1\|$ for any $w_1,w_2\in W$. Therefore, thanks to \eqref{S_via_hat_S}, we obtain the following result.

\begin{Proposition}\label{Prop1} Suppose that $F(\bar x,\bar w)=0$ and the Lagrange multiplier $\lambda$ corresponding to the stationary point $\bar x\in S(\bar w)$ is positive. If condition \eqref{Aubin_Suffi_Condi2} is satisfied, then $S$ has a Lipschitz continuous single-valued localization around $\bar w$ for $\bar x$. 
\end{Proposition}

Clearly, Proposition \ref{Prop1} encompasses Remark 4.1 of Part 1, which gives a sufficient condition for the local Lipschitz-like property of $S$ around $(\bar w,\bar x)$.

\medskip
{\sc Case 3}: $F(\bar x,\bar w)=0$ \textit{and} $\lambda =0$. In this case, using the results of Dontchev and Rockafellar \cite{DR_1996} one can verify the following lemma; see Sect.~\ref{appendices} for details. 
	 \begin{Lemma}\label{new_lemma2} The generalized equation \eqref{par_gen_eq} is strongly regular at $(\bar x,\lambda)$ iff the matrix $\nabla_{xx}^2{\mathcal{L}}(\bar x,\lambda,\bar w)=\nabla_{xx}^2f_0(\bar x,\bar w)$
is nonsingular and
\begin{equation}\label{strict_ineq(3)}\nabla_{x}F(\bar x,\bar w)^T\nabla_{xx}^2f_0(\bar x,\bar w)^{-1}\nabla_{x}F(\bar x,\bar w)>0.\end{equation}
\end{Lemma}

The sufficient condition for $S$ to be locally Lipschitz-like around $(\bar w,\bar x)$ in assertion (b) of Theorem 4.2 in Part 1 is \eqref{Aubin_Suffi_Condi3}, which 
reads as
\begin{equation}\label{4a}
{\rm ker}\,A'_1\cap {\rm ker}\,A'_2=\{0\} \tag{\text{4a}}
\end{equation}
\begin{equation}\label{4b}
{\rm ker}\,A'_1\cap ({\rm ker}\nabla_xF(\bar x,\bar w)\times \R)\subset {\rm ker}\,A'_2=\{0\} \tag{\text{4b}}
\end{equation}
\begin{equation}\label{4c}
{\rm ker}\,A'_1\cap\Delta_1\subset {\rm ker}\,A'_2\tag{\text{4c}}
\end{equation}
\begin{equation}\label{4d}
{\rm ker}\,\nabla_{xx}^2f_0(\bar x,\bar w)\cap\Delta_2\subset {\rm ker}\,\nabla_{wx}^2f_0(\bar x,\bar w) \tag{\text{4d}}
\end{equation}
where
$$A'_1=\left(\begin{matrix}
\nabla_{xx}^2f_0(\bar x,\bar w)&\ \, \nabla_xF(\bar x,\bar w)
\end{matrix}\right),$$
$$A'_2=\left(\begin{matrix}
\nabla_{wx}^2f_0(\bar x,\bar w)& \ \, \nabla_wF(\bar x,\bar w)
\end{matrix}\right),$$
$$\Delta_1=\{(v'_1,\gamma)\mid \nabla_{x}F(\bar x,\bar w)^Tv'_1>0,\;\gamma\geq 0\},$$
$$\Delta_2=\{(v'_1,\gamma)\mid \nabla_{x}F(\bar x,\bar w)^Tv'_1<0\}.$$
Now conditions \eqref{4a}--\eqref{4c} imply
\begin{equation}\label{5a}
{\rm ker}\,A'_1\cap ({\rm ker}\nabla_xF(\bar x,\bar w)\times \R)=\{0\},
\end{equation}
\begin{equation}\label{5b}
{\rm ker}\,A'_1\cap\Delta_1=\emptyset.
\end{equation} Indeed, by \eqref{4b} one sees that the linear subspace ${\rm ker}\,A'_1\cap ({\rm ker}\nabla_xF(\bar x,\bar w)\times \R)$ is contained in ${\rm ker}\,A'_1\cap {\rm ker}\,A'_2$. So, by \eqref{4a}, the subspace just consists of the origin. This justifies \eqref{5a}. Similarly, by \eqref{4c}, the set ${\rm ker}\,A'_1\cap\Delta_1$ is  contained in ${\rm ker}\,A'_1\cap {\rm ker}\,A'_2$. Hence, from  \eqref{4a} it follows that  ${\rm ker}\,A'_1\cap\Delta_1\subset\{0\}$. As $0\notin \Delta_1$, condition \eqref{5b} is satisfied.

Condition \eqref{5a} implies that $\det\nabla_{xx}^2f_0(\bar x,\bar w)\neq 0$. Indeed, if there existed $x'\in\R^n\setminus\{0\}$, then by choosing $\tau'=0$ we would have $$(x',\tau')\in {\rm ker}\,A'_1\cap ({\rm ker}\nabla_xF(\bar x,\bar w)\times \R).$$ This contradicts \eqref{5a}. 

To see that \eqref{5a} and \eqref{5b} yield \eqref{strict_ineq(3)}, put $v'=\nabla_{xx}^2f_0(\bar x,\bar w)^{-1}\nabla_{x}F(\bar x,\bar w))$. We have to show that $\nabla_{x}F(\bar x,\bar w)^Tv'>0$. If $\nabla_{x}F(\bar x,\bar w)^Tv'=0$, then 
$$(v',-1)\in {\rm ker}\,A'_1\cap ({\rm ker}\nabla_xF(\bar x,\bar w)\times \R).$$ This contradicts \eqref{5a}. If $\nabla_{x}F(\bar x,\bar w)^Tv'<0$, then for $v'_1:=-v'$ one has $\nabla_{x}F(\bar x,\bar w)^Tv_1'>0$. Choosing $\gamma=1$, by direct calculation we can verify that $(v'_1,\gamma)\in {\rm ker}\,A'_1\cap\Delta_1$. This is a contraction to \eqref{5b}. 

We have thus proved that the conditions $\det\nabla_{xx}^2f_0(\bar x,\bar w)\neq 0$ and \eqref{strict_ineq(3)} follow from \eqref{4a}--\eqref{4c}. Hence, if the conditions \eqref{4a}--\eqref{4c} are satisfied, then \eqref{par_gen_eq} is strongly regular at $(\bar x,\lambda)=(\bar x,0)$. Therefore,  invoking Theorem~2.1 from \cite{Robinson_1980} to the parametric generalized equation \eqref{par_gen_eq}, we can assert that if \eqref{4a}--\eqref{4c} are satisfied, then the implicit multifunction $\widehat S(.)$ has a Lipschitz continuous single-valued localization around $\bar w$ for $(\bar x,\lambda)=(\bar x,0)$. Thus, thanks to \eqref{S_via_hat_S}, we have the following result.

\begin{Proposition}\label{Prop2} Suppose that $F(\bar x,\bar w)=0$ and the Lagrange multiplier $\lambda$ corresponding to the stationary point $\bar x\in S(\bar w)$ is zero. If \eqref{4a}--\eqref{4c} are satisfied, then $S$ has a Lipschitz continuous single-valued localization around $\bar w$ for $\bar x$. 
\end{Proposition}

The result stated in Proposition~\ref{Prop2} is better than assertion (b) of Theorem~4.2 in Part 1, which says that if \eqref{Aubin_Suffi_Condi3} is fulfilled, i.e., \eqref{4a}--\eqref{4d} are valid, then $S$ is locally Lipschitz-like around $(\bar w,\bar x)$. 


\section{Appendices}\label{appendices}

{\it Proof of Lemma \ref{new_lemma1}}\ By the definition of Robinson \cite[p.~45]{Robinson_1980}, the strong regularity of \eqref{par_gen_eq} at $(\bar x,\lambda)$ is identical to the strong regularity of the affine variational inequality
\begin{equation}\label{linear_variational_ineq}
0\in A\left(\begin{matrix}
x\\\alpha
\end{matrix}\right)+\bar q+N_{\R^n\times\R_{+}}(x,\alpha),
\end{equation}
where
\begin{equation}\label{matrix_A} A:=\nabla_{(x,\alpha)}g(\bar x,\lambda,\bar w)=\begin{pmatrix} \nabla_{xx}^2{\mathcal{L}}(\bar x,\lambda,\bar w) & \ \; \nabla_{x}F(\bar x,\bar w)\\
-\nabla_{x}F(\bar x,\bar w)^T & \ \; 0
\end{pmatrix}\end{equation} and $$\bar q:=g(\bar x,\lambda,\bar w) - A\left(\begin{matrix}
\bar x\\\lambda\end{matrix}\right),$$  at $(\bar x,\lambda)$.

According to \cite[Theorem~1]{DR_1996}, the affine variational inequality \eqref{linear_variational_ineq} is strongly regular at $(\bar x,\lambda)$ if and only if 
the multifunction $L:\R^n\times\R\rightrightarrows \R^n\times\R$ with
$$L(q):=\left\{\left(\begin{matrix}
x\\\alpha
\end{matrix}\right)\, \mid \, 0\in A\left(\begin{matrix}
x\\ \alpha
\end{matrix}\right)+ q+N_{\R^n\times\R_{+}}(x,\alpha) \right\}$$ is locally Lipschitz-like around $(\bar q, (\bar x,\lambda))$. Furthermore, applying \cite[Theorem~2]{DR_1996}, we can assert that the latter is valid iff the\textit{ critical face condition} holds at $(\bar q, (\bar x,\lambda))$, i.e., \textit{for any  choice of closed faces $F_1$ and $F_2$ of the critical cone $K_0$ with $F_1\supset F_2$,
	\begin{equation}\label{critical_face_condi_inclusion}
	\big[u\in F_1-F_2,\ \, A^Tu \in (F_1-F_2)^*\big]\ \; \Longrightarrow\ \;  u=0.
	\end{equation}}
Here, $$K_0=K((\bar x,\lambda),v_0):=\big\{(x',\alpha')\in T_{\R^n\times\R_{+}}(\bar x,\lambda)\mid (x',\alpha')\perp v_0\big\},$$ 
with $$v_0:=-A\left(\begin{matrix}
\bar x\\\lambda
\end{matrix}\right)-\bar q\in N_{\R^n\times\R_{+}}(\bar x,\lambda).$$
Recall that a convex subset $F$ of a convex set $C\subset\R^p$ is a \textit{face} of $C$ if every closed line segment in $C$ with a relative interior point in $F$ has both endpoints in $F$. When $K_0$ is a linear subspace of $\R^n\times\R$, it has a unique closed face, namely itself. Then, the critical face condition is reduced to
\begin{equation}\label{reduced_critical_face_condi_inclusion}
\big[u\in K_0,\ \, A^Tu \perp K_0\big]\ \; \Longrightarrow\ \;  u=0.
\end{equation}
\hskip0.5cm In the case $\lambda>0$, the critical face is equivalent to the nonsingularity of the matrix in \eqref{combined_matrix}. Indeed,
the condition $\lambda>0$ implies $N_{\R^n\times\R_{+}}(\bar x,\lambda)=\{(0,0)\}$,  
$T_{\R^n\times\R_{+}}(\bar x,\lambda)=\R^{n}\times \R$, and $v_0=(0,0)$. It follows that $K_0=\R^{n}\times \R$. So, the critical face is reduced to \eqref{reduced_critical_face_condi_inclusion}, which is
$$A^Tu=0\ \; \Longrightarrow\ \; u=0.$$
The latter means that $A$ is nonsingular; or, equivalently, the matrix \eqref{combined_matrix} is nonsingular. 

Thus, we have proved that the generalized equation \eqref{par_gen_eq} is strongly regular at $(\bar x,\lambda)$ iff the matrix \eqref{combined_matrix} is nonsingular. $\hfill\Box$

\medskip	
\noindent {\it Proof of Lemma \ref{new_lemma2}.}\ The arguments described in the beginning of the proof of Lemma \ref{new_lemma1} show that  the generalized equation \eqref{par_gen_eq} is strongly regular at $(\bar x,\lambda)$ iff the critical face condition holds at $(\bar q, (\bar x,\lambda))$, i.e., for any  choice of closed faces $F_1$ and $F_2$ of the critical cone $K_0$ with $F_1\supset F_2$ the condition \eqref{critical_face_condi_inclusion} is fulfilled.

Since $\lambda=0$, $N_{\R^n\times\R_{+}}(\bar x,\lambda)=\{0\}\times\R_{-}$, and 
$$T_{\R^n\times\R_{+}}(\bar x,\lambda)=\R^{n}\times\R_{+}.$$
As $v_0\in N_{\R^n\times\R_{+}}(\bar x,\lambda),$ there are two situations: (a)
$v_0=(0,\beta)$ with $\beta<0$; (b)~$v_0=(0,0)$. If (a) occurs, then $K_0=\R^{n}\times\{0\}$. Since $K_0$ is a linear subspace, the critical face condition is reduced to \eqref{reduced_critical_face_condi_inclusion}. Using the formula for $A$ in \eqref{matrix_A}, one can easily show that \eqref{reduced_critical_face_condi_inclusion} is equivalent to the requirement that the matrix $\nabla_{xx}^2{\mathcal{L}}(\bar x,\lambda,\bar w)$ is nonsingular. As $\lambda=0$, one has $\nabla_{xx}^2{\mathcal{L}}(\bar x,\lambda,\bar w)=\nabla_{xx}^2f_0(\bar x,\bar w)$. So, \eqref{reduced_critical_face_condi_inclusion} is also equivalent to the condition saying that the matrix $\nabla_{xx}^2f_0(\bar x,\bar w)$ is nonsingular. Now, suppose that the situation (b) occurs. Then, $$K_0=K((\bar x,\lambda),v_0)=\R^{n}\times\R_{+}.$$  Obviously, $K_0$ has only two nonempty faces: $\R^n\times \{0\}$ and $\R^{n}\times\R_{+}$. 

For $F_1=F_2=\R^{n}\times\{0\},$ one has  $F_1-F_2=\R^{n}\times\{0\}$. Then, $$(F_1-F_2)^*=\{0\}\times\R$$ and \eqref{critical_face_condi_inclusion} is satisfied iff, for any $u'\in\R^n$, $$\nabla_{xx}^2{\mathcal{L}}(\bar x,\lambda,\bar w)u'=0\ \;\Longrightarrow\ \;u'=0.$$
As $\lambda=0$, it holds that $\nabla_{xx}^2{\mathcal{L}}(\bar x,\lambda,\bar w)=\nabla_{xx}^2f_0(\bar x,\bar w).$
Therefore, \eqref{critical_face_condi_inclusion} is valid iff, for any $u'\in\R^n$,
$$\nabla_{xx}^2f_0(\bar x,\bar w)u'=0\ \;\Longrightarrow\ \;u'=0.$$
This is equivalent to saying that $\nabla_{xx}^2f_0(\bar x,\bar w)$ is nonsingular. 

For $F_1=F_2=\R^{n}\times\R_{+}$, $F_1-F_2=\R^{n}\times\R$. Then, $(F_1-F_2)^*=\{0\}\times \{0\}$ and \eqref{critical_face_condi_inclusion} is satisfied iff the matrix \begin{equation*}A^T=\begin{pmatrix} \nabla_{xx}^2f_0(\bar x,\bar w) & \ \; -\nabla_{x}F(\bar x,\bar w)\\
\nabla_{x}F(\bar x,\bar w)^T & \ \; 0
\end{pmatrix}\end{equation*} is nonsingular, or, $A$ is nonsingular. 

For $F_1=\R^{n}\times\R_{+}$ and $F_2=\R^{n}\times\{0\}$, $F_1-F_2=\R^{n}\times\R_{+}$ and $(F_1-F_2)^*=\{0\}\times\R_{-}$. Then, \eqref{critical_face_condi_inclusion} is fulfilled iff
\begin{equation}\label{system_new}
\begin{cases}
\nabla_{xx}^2f_0(\bar x,\bar w)u'-\gamma\nabla_{x}F(\bar x,\bar w)=0\\
\nabla_{x}F(\bar x,\bar w)^Tu'\leq 0\\
u'\in\R^n,\ \gamma\geq 0
\end{cases} \ \;\Longrightarrow\ \;
\begin{cases}
u'=0\\
\gamma=0.
\end{cases}
\end{equation} 
\hskip0.5cm The proof of the ``necessity part" of Lemma \ref{new_lemma2} will be completed if we can show that \eqref{strict_ineq(3)} is valid. If \eqref{strict_ineq(3)} does not hold, then by putting
$$u'=\nabla_{xx}^2f_0(\bar x,\bar w)^{-1}\nabla_{x}F(\bar x,\bar w),$$
we have $$\nabla_xF(\bar x,\bar w)^Tu'=\nabla_xF(\bar x,\bar w)^T\nabla_{xx}^2f_0(\bar x,\bar w)^{-1}\nabla_{x}F(\bar x,\bar w)\leq 0.$$ So, for $\gamma=1$, one has
$$\begin{cases}
	\nabla_{xx}^2f_0(\bar x,\bar w)u'-\gamma\nabla_{x}F(\bar x,\bar w)=0\\
	\nabla_{x}F(\bar x,\bar w)^Tu'\leq 0\\
	u'\in\R^n,\ \gamma\geq 0.
\end{cases}$$
This contradicts \eqref{system_new}. We have thus proved \eqref{strict_ineq(3)} is valid.

To prove the ``sufficiency part" of Lemma \ref{new_lemma2}, we suppose that the matrix $\nabla_{xx}^2{\mathcal{L}}(\bar x,\lambda,\bar w)=\nabla_{xx}^2f_0(\bar x,\bar w)$
 is nonsingular and \eqref{strict_ineq(3)} is fulfilled. To verify the fulfillment of the critical face condition at $(\bar q, (\bar x,\lambda))$, we need only to show that the matrix $A$ is nonsingular and the implication \eqref{system_new} holds. 
 
 To obtain the nonsingularity of $A$, suppose to the contrary that there exists a pair $(u',\gamma)\neq (0,0)$ satisfying
 \begin{equation}\label{two_equalities}
 \begin{cases}
 \nabla_{xx}^2f_0(\bar x,\bar w)u'-\gamma\nabla_{x}F(\bar x,\bar w)=0\\
 \nabla_{x}F(\bar x,\bar w)^Tu'=0.
 \end{cases}
 \end{equation}
If $\gamma=0$, then the first equation of \eqref{two_equalities} forces $u'=0$, because the matrix $\nabla_{xx}^2f_0(\bar x,\bar w)$ is nonsingular. So, we must have $\gamma\neq 0$. From the first equation of \eqref{two_equalities}  we deduce that
$$u'=\gamma\nabla_{xx}^2f_0(\bar x,\bar w)^{-1}\nabla_{x}F(\bar x,\bar w).$$ Hence, by the second equation of \eqref{two_equalities}, we obtain
 $$\gamma\nabla_{x}F(\bar x,\bar w)^T\nabla_{xx}^2f_0(\bar x,\bar w)^{-1}\nabla_{x}F(\bar x,\bar w)=0.$$ This obviously contradicts \eqref{strict_ineq(3)}. Thus, $A$ is nonsingular.
 
 Finally, to obtain the implication \eqref{system_new},  let $u'\in\R^n$ and $\gamma\geq 0$ be such that 
 \begin{equation}\label{system_new_1}
 \begin{cases}
 \nabla_{xx}^2f_0(\bar x,\bar w)u'-\gamma\nabla_{x}F(\bar x,\bar w)=0\\
 \nabla_{x}F(\bar x,\bar w)^Tu'\leq 0
 \end{cases} 
 \end{equation} 
 Multiplying both sides of the equation in \eqref{system_new_1} from the left with the $1\times n$ matrix $\nabla_{x}F(\bar x,\bar w)^T\nabla_{xx}^2f_0(\bar x,\bar w)^{-1}$, one obtains
\begin{equation}\label{special_equation}\nabla_{x}F(\bar x,\bar w)^Tu'-\gamma\nabla_{x}F(\bar x,\bar w)^T\nabla_{xx}^2f_0(\bar x,\bar w)^{-1}\nabla_{x}F(\bar x,\bar w)=0.\end{equation}
Due to \eqref{strict_ineq(3)} and the condition $\gamma\geq 0$, 
$$-\gamma\nabla_{x}F(\bar x,\bar w)^T\nabla_{xx}^2f_0(\bar x,\bar w)^{-1}\nabla_{x}F(\bar x,\bar w)\leq 0.$$
Combining this with the inequality $\nabla_{x}F(\bar x,\bar w)^Tu'\leq 0$ from \eqref{system_new_1}, by \eqref{special_equation} one has 
$$
-\gamma\nabla_{x}F(\bar x,\bar w)^T\nabla_{xx}^2f_0(\bar x,\bar w)^{-1}\nabla_{x}F(\bar x,\bar w)=0.$$
Due to \eqref{strict_ineq(3)}, $\gamma=0$. Then, the first equation in \eqref{system_new_1} implies $\nabla_{xx}^2f_0(\bar x,\bar w)u'=0$. As $\nabla_{xx}^2f_0(\bar x,\bar w)$ is nonsingular, one has $u'=0$. Thus, \eqref{system_new} is valid. 

The proof is complete. $\hfill\Box$

\section{Conclusions}

In this paper, we have analyzed the stability of the stationary point set map of a smooth parametric optimization problem with one smooth functional constraint under total perturbations. Not only sufficient conditions for the local Lipschitz-like property of the stationary point set map were given, but also necessary conditions for the latter property have been obtained. Sufficient conditions for the Robinson stability of the stationary point set map were given. In addition, we have revisited several stability theorems in indefinite quadratic programming.

It is still unclear to us whether the Lipschitz continuous single-valued localization mentioned in Propositions \ref{Prop1} and \ref{Prop2} implies the Robinson stability of the stationary point set map, or not.

Extensions of the obtained results to the case of smooth parametric optimization problem with more than one smooth functional constraint under total perturbations are worthy further investigations.

 \begin{acknowledgements}	
 	This work was supported by National Foundation for Science $\&$ Technology Development (Vietnam) and the Grant MOST 105-2115-M-039-002-MY3 (Taiwan). The authors are grateful to the anonymous referees for their careful readings and valuable suggestions. Section 5 is based on the comments made by one of the referees.
 \end{acknowledgements}

\end{document}